\documentclass[12pt]{article}

\usepackage{amsmath}
\usepackage{amssymb}
\usepackage{dsfont}
\usepackage{natbib}
\usepackage{abstract}
\usepackage{hyperref}

\newcommand{\beqn}{\begin{eqnarray*}}
\newcommand{\eeqn}{\end{eqnarray*}}
\newcommand{\bneqn}{\vspace{-0.25cm}\begin{eqnarray}}
\newcommand{\eneqn}{\end{eqnarray}}
\newcommand{\parens}[1]{\left(#1\right)}
\newcommand{\squared}[1]{\parens{#1}^2}
\newcommand{\bracks}[1]{\left[#1\right]}

\newcommand{\expe}[1]{\mathbb{E}\bracks{#1}}
\newcommand{\cexpe}[2]{\expe{#1 ~ | ~ #2}}

\newcommand{\reals}{\mathbb{R}}

\title{The General Stationary Gaussian Markov Process}

\author{Larry Brown,$^1$ Philip Ernst,$^1$ Larry Shepp$^1$ and Robert Wolpert$^2$ \\~\\ $^1$Department of Statistics, \\ The Wharton School of the University of Pennsylvania \\ $^2$Department of Statistical Science, Duke University}

\begin{document}
\maketitle

\begin{abstract}
We find the class, ${\cal{C}}_k, k \ge 0$, of all zero
mean stationary Gaussian processes, $Y(t), ~t \in \reals$ with $k$
derivatives, for which 

\beqn
Z(t) \equiv (Y^{(0)}(t), Y^{(1)}(t), \ldots, Y^{(k)}(t) ), ~ t \ge 0
\eeqn

\noindent
is a $(k+1)$-vector Markov process. (here, $Y^{(0)}(t) = Y(t)$).
\end{abstract}

\section{Introduction}
\noindent
We show that the process, $Y$, can be described in three equivalent ways:\\

\noindent
(i). Each member of ${\cal{C}}_k$ is given uniquely by a certain polynomial
$P(z)$ via the covariance of any such $Y$:

\beqn
\expe{Y(s)Y(t)} = R(s,t) = r(t-s) = \int_{-\infty}^\infty \frac{e^{i(t-s)z dz}}{|P(z)|^2}
\eeqn

\noindent
where $P(z)$ is a polynomial of degree $k+1$ with positive leading coefficiant
and all complex roots $\zeta_j = \rho_j + i \sigma_j,j=0,\ldots,k$ with
$\sigma_j > 0$, and $\rho_j$ real. If some $\rho_j \ne 0$, then there is
another $\zeta_k = -\zeta_j^*$ which is the negative conjugate of $\zeta_j$.
Note $r(t)$ automatically has $2k$ derivatives at each $t$, but not $2k+1$.\\

\noindent
(ii) Equivalently, it is necessary and sufficient that $Y \in {\cal{C}}$ has
the representions via Wiener integrals with a standard Brownian motions, $W$

\beqn
Y(t) = \int_{-\infty}^\infty g(t-\theta) dW(\theta)
\eeqn

\noindent
where $g$ has $L^2$ Fourier transform, 
\beqn
{\hat{g}}(z) = \frac{1}{|P(z)|}
\eeqn
or,
equivalently, $Y$ has the spectral representation, with a pair of independent
standard Brownian motions, $W_1,W_2$, with $f(z) \equiv {\hat{g}}$,

\beqn
Y(t) = \int_{-\infty}^\infty \cos{tz} f(z) dW_1(z) + \int_{-\infty}^\infty \sin{tz} f(z) dW_2(z)
\eeqn

\noindent
(iii) Equivalently, it is necessary and sufficient for $Y \in {\cal{C}}$ that
$Z$ has a representation as an Ito vector diffusion process:

\beqn
dY^{(i)}(t) &=& Y^{(i+1)}(t) dt,~ 0 \le i < k, t \in \reals \\
dY^{(k)}(t) &=& \sum_{j=0}^k a_j Y^{(j)}(t) dt + b dW(t), ~t \in \reals
\eeqn

\noindent
where the coefficients, $a_j$, are the unique solution of the equations:

\beqn
r^{(k+i+1)}(0^+) = \sum_{j=0}^k a_j r^{(i+j)}(0), ~i=0, 1, \ldots, k
\eeqn

\noindent
Note that the left and right derivatives of $r^{(j)}$ are equal except for
$j = 2k$. The diffusion coefficient, $b$, is given by

\beqn
b^2 = \sum_{j=0}^k a_j r^{(j+k)}(0)(-1)^{j+1} + (-1)^k r^{(2k+1)}(0^-)
\eeqn

\noindent
Given the polynomial, $P(z)$, in (i), which uniquely determines each process
$Y \in {\cal{C}}$, the derivatives $r^{(j)}(0)$ are easily obtained from the
representation in (i) above as

\beqn
r^{(j)}(0) = \int_{-\infty}^\infty \frac{(iz)^j dz}{|P(z)|^2}, ~0 \le j \le 2k+1
\eeqn

\noindent
where for $j = 2k+1$ the integral is not $L^1$ convergent, but is understood
as a principal value. Then the coefficients $a_j, b$ of the Ito equation is
determined as indicated above. This allows one to determine exactly which
Ito vector equations have stationary solutions, and what the stationary
distribution is. Of course it is the Gaussian vector, $Y^{(j)}(0)$, with
covariance 

\beqn
r^{(i+j)}(0), ~i,j = 0,\ldots,k
\eeqn

\noindent
It seems that these results are new despite the fact that the problem has
been around for nearly fifty years to describe the class ${\cal{C}}$.\\

\section{The case $k = 0$}.

For $k = 0$, we ask what is the set of all stationary Gauss-Markov
processes and the answer is well-known as the set of processes with covariance
$r(t) = Ae^{-\alpha |t|},~ A \ge 0,~ \alpha \ge 0$, with representation

\beqn
Y(t) = \sqrt{A} e^{-\alpha t}W\parens{e^{2\alpha t}}
\eeqn

\noindent with $W$ a standard Wiener
process, which satisfies an Ito equation of the form

\beqn
dY(t) = a Y(t) dt  + b dW(t), ~~\text{where}~~ a < 0 ~~ \text{and} ~~ b > 0
\eeqn

Let us give the simple proof for $k = 0$ that the covariance $r$ must be as
stated, because we will use the same method for the general case, $k$, and this
will make things clearer. The idea is to find or define the covariance, 
$R(s,t) = r(|t-s|)$, under the stated assumptions. We must have $Y(s)$ and
$Y(u)$ given $Y(t)$ to be uncorrelated for Markovianness to hold, whenever
$s < t < u$. Since 

\beqn
\cexpe{Y(u)}{Y(t)} = Y(t) \frac{r(t-s)}{r(0)}
\eeqn

\noindent
this means that $r$ satisfies for any positive $u,v$,

\beqn
\expe{Y(u) - \frac{r(u)}{r(0)} Y(0)} Y(-v) = 0, ~~\text{or} ~~ r(u+v) = \frac{r(u) r(v)}{r(0)}
\eeqn

\noindent
Since $r$ is continuous and nonnegative definite, it follows that
$r(h) = Ae^{-h \alpha},$ $h \ge 0$. To see this note that 
$f(t) = \log{\frac{r(t)}{r(0)}}$ satisfies $f(u+v) = f(u) + f(v)$ and so
$f(\frac{m}{n}) = f(1)\frac{m}{n}$. Since $f$ is continuous we have 
$f(u) = u f(1)$ and the claim follows. Since $r$ is to be nonnegative definite
we must have $A \ge 0$ and $\alpha \ge 0$. We have found a necessary condition
for $r$ to be the covariance. In fact we see that $r$ is infinitely 
differentiable, except at $u = 0$, where $r$ has finite left and right
derivatives. An analogous property will hold for every $k$. The only thing
missing is to prove sufficiency, i.e., the existence of a process with this
covariance. This is easy if we use the representation:

\beqn
Y(t) = \sqrt{A}e^{-t \alpha}W(e^{2\alpha t})
\eeqn

\noindent which has covariance $r(t-s) = A e^{-|t-s|\alpha}$.\\

\noindent
We know $Y$ is Markovian, so we can obtain the coefficients of the Ito equation
by the formula,

\beqn
aY(0) = \lim_{h \downarrow 0} \frac{\cexpe{Y(h) - Y(0)}{Y(0)}}{h} = \lim_{h \downarrow 0} \frac{Y(0) (r(h)-r(0))}{h r(0)} = -\alpha Y(0)
\eeqn

\noindent
and the other Ito coefficient is

\beqn
b = \lim_{h \downarrow 0} \frac{\squared{\expe{Y(h}- \cexpe{Y(h)}{Y(0)}}}{h} = \lim_{h \downarrow 0} \frac{r(0) - \frac{r^2(h)}{r(0)} }{h} = 2 \alpha A
\eeqn

\noindent 
so the Ito equation is 

\beqn
dY(t) = -\alpha Y(t) dt + \sqrt{2A \alpha} dW(t)
\eeqn

\noindent
Finally, the stationary measure has $\sigma^2 = r(0) = A$.

\section{The case $k = 1$}.

We will give the approach for general $k$ but let's do $k = 1$. Since $Y(u)$
and $Y(s)$ are conditionally uncorrelated given $Y(t)$ for $s < t < u$, we need

\beqn
(\expe{Y(u)} - \cexpe{Y(u)}{Y(0),Y^\prime(0)})Y(-v) = 0,
\eeqn

\noindent
for $u > 0 > -v$, and since we must have $r^\prime(0) = 0$ since $r$ is even
and is differentiable at zero because $Y$ is differentiable, we have

\beqn
\cexpe{Y(u)}{Y(0),Y^\prime(0)} = \frac{r(u))}{r(0)} Y(0) + \frac{r^\prime(u)}{r^{\prime \prime}(0)} Y^\prime(0)
\eeqn

\noindent
Since $(\expe{Y(u)} - \cexpe{Y(u)}{Y(0),Y^\prime(0)})Y(-v) = 0$, for $u > 0, v > 0$, this
gives

\beqn(*) \hskip 1in r(u+v) -\frac{r(u) r(v)}{r(0)} - \frac{r^\prime(u) r^\prime(v)}{r^{\prime \prime}(0) } = 0, u > 0, v > 0
\eeqn

\noindent
This shows that $r(u), u > 0$ is infinitely differentiable because $r$ is
differentiable and $(*)$ exhibits $r^\prime$ in terms of $r$, so that 
$r^\prime$ is differentiable, and by induction $r$ has all derivatives at
$u \ne 0$. Moreover, since $r$ is even we must have $r^\prime(0) = 0$. We next
show that $r(u)$ satisfies a second order ODE with constant coefficients,
namely:

\beqn
r^{(2)}(u) = r^{(0)}(u) \frac{r^{(2)}(0)}{r^{(0)}(0)} + r^{(1)}(u) \frac{r^{(3)}(0+)}{r^{(2)}(0)}
\eeqn

\noindent
We claim first that $r$ is twice differentiable at $u = 0$. This follows from
the fact that $Y^\prime(t)$ is a Gaussian variable and so

\beqn
r^{\prime \prime}(0) = -\expe{Y^\prime(0) Y^\prime(0)}
\eeqn

exists.

\noindent
If we expand each side of $(*)$ into power series in $v$, then we get that
up to a term, $o(v^2)$,

\beqn
\sum_{j=0}^2 \frac{r^{(j)}(u) v^j}{j!} = \frac{r(u)}{r(0)} \sum_{j=0}^2 \frac{r^{(j)}(0) v^j}{j!} + \frac{r^{(1)}(u)}{r^{(2)}(0)} \sum_{j=0}^2 \frac{r^{(j+1)}(0) v^j}{j!}
\eeqn

\noindent
and we see that the coefficients of $v^0,v^1$ vanish automatically, but the
coefficient of $v^2$ shows that $r^{(2)}(0+)$ exists and then the coefficient 
$v^2$ being equal on both sides gives a second degree differential equation
for $r$, for $u > 0$. Thinking of $r(0), r^{(2)}(0),r^{(3)}(0+)$ as constants,
we see that for $u > 0$, the differential equation for $r$ has constant
coefficients. If the indicial equation has distict roots, then this means that

\begin{equation}
r(u) = \sum_{j=1}^2 A_j e^{-u a_j}. \label{eqn:rep1}
\end{equation}

If the two roots are not distinct, then one gets a limiting covariance e.g.,
for the case when $a_1 = a, a_2 = a+\epsilon$, and
$A_1 = -A_2 = \frac{1}{\epsilon}$, $r$ becomes the derivative,

\beqn
r(u) = (1+au)e^{-ua}, u \ge 0
\eeqn

\noindent
Next, we have to check that any such $r$ satisfies equation (*). This is
easy to check in this case, so that satisfying (*) imposes no additional
restrictions than satisfying a second order ode with constant coefficients.
It is not true that every such $r$ is realizable because $r$ must be a
covariance and conditions to ensure this must be placed on $a_j$ and $A_j$.
For example we must have $a_j > 0$. For $Y$ to be differentiable we
need that $r(h)$ be twice differentiable at $h = 0$. In turn, this means that
$-r^\prime(0) = A_1 a_1 + A_2 a_2 = 0$. We also need that $r(0) = A_1+A_2 > 0$.
Let us use the process represented below to show that with these restrictions
the condition is also sufficient to realize the covariance in (\ref{eqn:rep1}):

\beqn
Y(t) = \int_{-\infty}^\infty f(t-\theta) dW(\theta)
\eeqn

\noindent
where $f$ is any $L^2$ function to get a class of processes with covariance
$r$ of the form above. Set:

\beqn
f(x) = A^- e^{x a^-}, ~x < 0, ~f(x) = A^+ e^{-x a^+}, ~x > 0
\eeqn

\noindent
Now, $Y$ will only be well defined when $f \in L^2$, so we need $a^\pm > 0$.
The covariance of the representation is easily seen to be

\beqn
r(u) = \parens{\frac{(A^-)^2}{2 a^-} - \frac{A^- A^+}{a^- - a^+}} e^{-h a^-}
  + \parens{\frac{(A^+)^2}{2 a^+} + \frac{A^- A^+}{a^- - a^+}} e^{-h a^+}
\eeqn

Also $Y$ will only be differentiable when $f$ is continuous. so this means 
$f(0-) = f(0+)$ or $A^- = A^+$. We get a certain class of covariances of our
form. Wolog we can choose $a_1 = a^-, a_2 = a^+$. If we set $A^+ = A^- = A$,
we need to choose $A$ so that

\beqn
A_1 = A^2\parens{\frac{1}{2a_-} - \frac{1}{a_- - a_+}},~ A_2 = A^2 \parens{\frac{1}{2a_+} + \frac{1}{a_- - a_+}}
\eeqn

\noindent
It is easy to check that $a_1 A_1 + a_2 A_2 = 0$ holds.\\

\subsection{Discrete Considerations}
Using definition (iii) to classify the elements ${\cal{C}}_k$, we ask for the form of an AR(2) process that will give rise to the continuous AR(2) process.

\section{General $k \ge 2$}.\\

We must have $Y(s)$ and $Y(u)$ conditionally independent given
$Y^{(j)}(t), j=0,\ldots, k$, and since the process is Gaussian this means
that $Y(s)$ and $Y(t)$ are conditionally uncorrelated. This means that

\beqn
\expe{Y(u) - \sum_{j=0}^k \alpha_j(u) Y^{j)}(0)} Y(-v) = 0
\eeqn

\noindent
for $u > 0, v > 0$ where

\beqn
\cexpe{Y(u)}{Y^{(j)}(0),j=0,\ldots,k} = \sum_{j=0}^k \alpha_j(u) Y^{(j)}(0)
\eeqn

\noindent
because conditional expectations are linear for Gaussian processes. Note
the $\alpha_j$'s are defined uniquely by the equations

\beqn
r^{(i)}(u) = \sum_{j=0}^k \alpha_j(u) r^{(i+j)}(0),~ i = 0,\ldots, k. \label{eqn:alpha}
\eeqn

\noindent
We would like to show that, without any further assumptions than the fact that
$r$ satisfies an ode of degree $k+1$, that the first equation holds, because
then we can conclude that the first equation poses no additional restrictions.
The first equation is the same as

\beqn
r(u+v) = \sum_{j=0}^k \alpha_j(u) r^{(j)}(v), ~u > 0, ~v > 0 \label{eqn:rank}
\eeqn

\noindent
i.e., $r(u+v)$ is of rank $k+1$, i.e., (\ref{eqn:rank}) holds if the
$\alpha_j$'s are defined by (\ref{eqn:alpha}). Here is where using the first
approach pays off to avoid a lot of algebra. Imagine solving (\ref{eqn:alpha})
for the $\alpha_j(u)$'s and then placing these $\alpha_j(u)$ into 
(\ref{eqn:rank}). We now let $v$ be small and positive and use power series.
We get that $r(u)$ satisfies a differential equation of degree $k+1$ with
constant coefficients as in the cases $k = 0,1$. But to avoid checking that
no further condition is required to prove that $r$ also satisfies 
(\ref{eqn:rank}) we can argue as follows. Let for $v$ fixed,

\beqn
f(u) = r(u+v) - \sum_{j=0}^k \alpha_j(u) r^{(j)}(v), ~u > 0
\eeqn

Note that $r(u+v)$ and $r^{(j)}(u)$, as functions of $u$, for all $j$ and any
fixed $v$ satisfy the same differential equation because the differential
equation has constant
coefficients. Also, there are $k+1$ zero boundary conditions $f^{(j)}(0+) = 0$,
so that $f \equiv 0$. This is quite subtle, and we need that $u > 0$ and
$v > 0$ here to get the required differentiability. It follows that
(\ref{eqn:rank}) holds.

So we have proved that if $r$ is a covariance for which $Z$ is a 
$(k+1)$-vector Markov process, then $r(u), u > 0$ satisfies a differential
equation of degree $k+1$ with constant coefficients. The general solution
for $r(u)$ must also be of the form (since $r$ is nonengative definite):

\beqn
r(u) = \int_{\reals} e^{iux} \mu(dx)
\eeqn

\noindent
for some nonnegative (spectral) measure, $\mu$. Since $r(u) = r(-u)$, $\mu$ 
is even, and since $r$ satisfies a differential equation of order $k+1$ with
constant coefficients, we must have for $u \ne 0$,

\beqn
\sum_{j=0}^k b_j r^{(j)}(u) = 0 = \int_{\reals} \sum_{j=0}^k b_j (-ix)^j e^{iux} \mu(dx) = \int_{\reals} e^{ixu} P(x) \mu(dx)
\eeqn

\noindent
We next prove that we must have, with $c > 0$,

\beqn
r(t) =\int_{\reals} \frac{e^{izt}}{|P(z)|^2} dz
\eeqn

\noindent
where $P(z) = c \prod_{j=0}^k (1- \frac{z}{\zeta_j})$, with $c > 0$, and with
the $k+1$ complex numbers, $\zeta_j, j = 0,\ldots,k$ having strictly positive
imaginary part. We require that $f(z) = f(-z)$ so we must have for each
$\zeta_j$ another $\zeta_{j^\prime}$ for which $\zeta_j = -\zeta_{j^\prime}^*$
is the negative complex conjugate. It may be that $j^\prime = j$ in which case
$\zeta_j$ is on the positive imaginary axis. For such a polynomial $P$, there
is an ode with constant coefficients satisfied by $r(t),t > 0$ because, by
Cauchy's theorem, the differential operator $P(-iD) r(t)$, $D = \frac{d}{dt}$,
is just

\beqn
P(-iD)r(t) = \int_{\reals} \frac{P(z)}{P(z)P^*(z)}e^{izt} dz =\int_{\reals} \frac{e^{itz}}{P^*(z)} dz = 0
\eeqn

\noindent
because we can complete the integral by adding a semicircle above the real 
$z$-axis along which, for $u > 0$, $e^{iuz}$ is bounded, and since the last
integrand is analytic in the upper half plane the integral is zero, and as the
radius of the semicircle goes to infinity the contribution from the arc is 
negligible because $P^*(z)$ is large.

The representation of $r$ as a covariance is now immediate, since we can just
set

\beqn
Y(t,\omega) = \int_{\reals} \cos{tz} \frac{dW_1(z,\omega)}{|P(z)|} + \int_{\reals} \sin{tz} \frac{dW_2(z,\omega)}{|P(z)|}
\eeqn

\noindent
where $W_i$ are iid standard Brownian motions, and check that $Y$ has
covariance $r$. Note that we cannot have any polynomial factor in the
numerator of the above equation for $r$ because $r^{(2k)}(u)$ must exist.

The Ito equation for 
\beqn
Z(t) = (Y^{(0)}(t),Y^{(1)}(t),\ldots,Y^{(k)}(t) )
\eeqn
is degenerate for the coefficients of $dY^{(j)}(t), j < k$, since

\beqn
dY^{(i)}(t) \equiv Y^{(i+1)}(t) dt, i < k
\eeqn

\noindent
but for $i = k$, we need to compute $a_j$ and $b$ in

\beqn
dY^{(k)}(t) = \sum_{j=0}^k a_j Y^{(j)}(t) dt + b dW(t)
\eeqn

\noindent
The $a_j$'s are found as follows, we may as well take $t = 0$, so

\beqn
\sum_{j=0}^k a_j Y^{(j)}(0) &=& \lim_{h \downarrow 0} \frac{\cexpe{Y^{(k)}(h)-Y^{(k)}(0)}{Y^{(0)}(0),\ldots,Y^{(k)}(0)}}{h}\\
&=&\lim_{h \downarrow 0} \frac{\sum_{j=0}^k \alpha_j^{(k)}(h)Y^{(j)}(0)-Y^{(k)}(0)}{h} = \sum_{j=0}^k \alpha_j^{(k+1)}(0)Y^{(j)}(0)
\eeqn

\noindent
where we have used the fact that $\alpha_j^{(i)}(0) = \delta_{i,j}$ because 
$\alpha_j^{(i)}(h)$ satisfies, for any $i \ge 0$,

\beqn
r^{(i)}(h) = \sum_{j=0}^k\alpha_j^{(i)}(h)r^{(j)}(0)
\eeqn

\noindent
and we may set $h = 0$. Comparing coefficeints, we can read off the result,

$a_j = \alpha_j^{(k+1)}(0)$.

\noindent
Moreover, the values of $a_j$ are given as the solutions of the equations
satisfied by the $\alpha_j^{(i)}$'s, i.e.,

\beqn
r^{(i+k+1)}(0) = \sum_{j=0}^k a_j r^{(i+j)}(0), i = 0,1,\ldots,k
\eeqn

\noindent
so the values of $a_j$ are now determined. To determine $b$, we use

\beqn
b^2 &=& \lim_{h \downarrow 0} \frac{\squared{\expe{Y^{(k)}(h)} - \cexpe{Y^{(k)}(h)}{Y^{(0)},\ldots,Y^{(k)}(0)}}}{h}\\ &=& \lim_{h \downarrow 0} \frac{r^{(2k)}(0)(-1)^k -\sum_{j=0}^k \alpha_j^{(k)}(h)r^{(k+j)}(h)(-1)^j}{h}\\
&=& \sum_{j=0}^{k-1} \alpha_j^{(k+1)}(0) r^{(k+j)}(0)(-1)^{j+1} + (-1)^{k+1}\lim_{h \downarrow 0} \frac{r^{(2k)}(h) \alpha_k^{(k)}(h) - r^{(2k)}(0)\alpha_k^{(k)}(0)}{h}.
\eeqn

\noindent
but the last limit is just the value of the derivative of the product of
$r^{(2k)}(h)\alpha_k^{(k)}(h)$ at $h = 0$, so we finally arrive at

\beqn
b^2 = \sum_{j=0}^k \alpha_j^{(k+1)}(0)r^{(k+j)}(0)(-1)^{j+1}+(-1)^{k+1}r^{(2k+1)}(0^+)
\eeqn

\noindent
or in terms of the already calculated $a_j$'s, switching to $0^-$

\beqn
b^2 = \sum_{j=0}^k a_j r^{(k+j)}(0)(-1)^{j+1} + (-1)^k r^{(2k+1)}(0^-)
\eeqn.

\noindent
All the coefficients can be computed in terms of the unique polynomial $P$
which corresponds to any process $Y$ in ${\cal{C}}$ because the only thing we
need to know are the derivatives of $r(t)$ at $t =0$, which are given by

\beqn
r^{(j)}(0^\pm) = \mp \int_{\reals} \frac{(iz)^j dz}{|P(z)|^2}
\eeqn

\noindent
where the integral is absolutely convergent for $j \le 2k$ and is understood 
as a principle value integral for $j = 2k+1$.\\

\noindent \textbf{Generalizations}

The problem also makes sense for $k = \infty$: the paths of $Y$ are then entire
analytic functions. In the case $k = \infty$, the Markov process degenerates 
because the values of $Y^{(j)}(0)$ for all $j$ completely determines the past
and the future of $Y$ because of the power series representation,

\beqn
Y(t) = \sum_{j=0}^\infty \frac{Y^{(j)}(0)t^j}{j!}
\eeqn.

\noindent Another generalization is to allow $Y$ itself to be a vector process, 

\beqn
{\bf{Y}}(t) = (Y_1(t),\ldots,Y_n(t) )
\eeqn

\noindent
and ask the same question. When is ${\bf{Y}}$ together with its first $k$
derivatives a mean zero stationary Gaussian Markov process. It seems there is
no real trouble making this generalization, although noncommuting matrices
enter. Indeed, the equation for the covariance,

\beqn
{\bf{R}}(s,t) = \expe{{\bf{Y}}(s){\bf{Y}}(t)} = {\bf{r}}(t-s)
\eeqn

even when $k = 0$ is

\beqn
r_{ij}(u+v) = \sum_{k,l} r_{ik}(u)A_{kl}r_{lj}(v)
\eeqn

\noindent
for some matrix $A$ which may not commute with the matrix $r_{ij}(u)$.

\end{document}